\documentclass[12pt]{amsart}
\usepackage{tikz-cd}
\usepackage{amssymb,amsmath}
\usepackage{array}

\usepackage{amsmath, amssymb, amsfonts, amstext, amsthm, amscd}
\usepackage[mathscr]{euscript}
\usepackage{enumerate}

\usepackage{float}
\usepackage{graphicx}

\usepackage[raiselinks=false,colorlinks=true,citecolor=blue,urlcolor=blue,
linkcolor=blue,bookmarksopen=true,pdftex]{hyperref}

\usepackage{tikz}
\usepackage{tikz-cd}
\usetikzlibrary{snakes}
\usetikzlibrary{intersections, calc}

\usepackage[english]{babel}
\usepackage{chngcntr}

\begin{document}

\thispagestyle{empty}

\def\theequation{\arabic{section}.\arabic{equation}}

\newcommand{\codim}{\mbox{{\rm codim}$\,$}}
\newcommand{\stab}{\mbox{{\rm stab}$\,$}}
\newcommand{\lr}{\mbox{$\longrightarrow$}}

\newcommand{\be}{\begin{equation}}
\newcommand{\ee}{\end{equation}}

\newtheorem{guess}{Theorem}[section]
\newcommand{\bth}{\begin{guess}$\!\!\!${\bf }~}
\newcommand{\eeth}{\end{guess}}
\renewcommand{\bar}{\overline}
\newtheorem{propo}[guess]{Proposition}
\newcommand{\bpropo}{\begin{propo}$\!\!\!${\bf }~}
\newcommand{\epropo}{\end{propo}}

\newtheorem{lema}[guess]{Lemma}
\newcommand{\blem}{\begin{lema}$\!\!\!${\bf }~}
\newcommand{\elem}{\end{lema}}

\newtheorem{defe}[guess]{Definition}
\newcommand{\bdefe}{\begin{defe}$\!\!\!${\bf }~}
\newcommand{\edefe}{\end{defe}}

\newtheorem{coro}[guess]{Corollary}
\newcommand{\bcor}{\begin{coro}$\!\!\!${\bf }~}
\newcommand{\ecor}{\end{coro}}

\newtheorem{rema}[guess]{Remark}
\newcommand{\brem}{\begin{rema}$\!\!\!${\bf }~\rm}
\newcommand{\erem}{\end{rema}}

\newtheorem{exam}[guess]{Example}
\newcommand{\beg}{\begin{exam}$\!\!\!${\bf }~\rm}
\newcommand{\eeg}{\end{exam}}

\newtheorem{notn}[guess]{Notation}
\newcommand{\bnot}{\begin{notn}$\!\!\!${\bf }~\rm}
\newcommand{\enot}{\end{notn}}

\newcommand{\ch}{{\mathcal H}}
\newcommand{\cf}{{\mathcal F}}
\newcommand{\cd}{{\mathcal D}}
\newcommand{\cR}{{\mathcal R}}
\newcommand{\cv}{{\mathcal V}}
\newcommand{\cn}{{\mathcal N}}
\newcommand{\lra}{\longrightarrow}
\newcommand{\ra}{\rightarrow}
\newcommand{\blr}{\Big \longrightarrow}
\newcommand{\da}{\Big \downarrow}
\newcommand{\ua}{\Big \uparrow}
\newcommand{\hra}{\mbox{{$\hookrightarrow$}}}
\newcommand{\rt}{\mbox{\Large{$\rightarrowtail$}}}
\newcommand{\dua}{\begin{array}[t]{c}
\Big\uparrow \\ [-4mm]
\scriptscriptstyle \wedge \end{array}}
\newcommand{\ctext}[1]{\makebox(0,0){#1}}
\setlength{\unitlength}{0.1mm}
\newcommand{\cl}{{\mathcal L}}
\newcommand{\cp}{{\mathcal P}}
\newcommand{\ci}{{\mathcal I}}
\newcommand{\bz}{\mathbb{Z}}
\newcommand{\cs}{{\mathcal s}}
\newcommand{\ce}{{\mathcal E}}
\newcommand{\ck}{{\mathcal K}}
\newcommand{\cz}{{\mathcal Z}}
\newcommand{\cg}{{\mathcal G}}
\newcommand{\cj}{{\mathcal J}}
\newcommand{\cc}{{\mathcal C}}
\newcommand{\ca}{{\mathcal A}}
\newcommand{\cb}{{\mathcal B}}
\newcommand{\cx}{{\mathcal X}}
\newcommand{\co}{{\mathcal O}}
\newcommand{\bq}{\mathbb{Q}}
\newcommand{\bt}{\mathbb{T}}
\newcommand{\bh}{\mathbb{H}}
\newcommand{\br}{\mathbb{R}}
\newcommand{\bl}{\mathbf{L}}
\newcommand{\wt}{\widetilde}
\newcommand{\im}{{\rm Im}\,}
\newcommand{\bc}{\mathbb{C}}
\newcommand{\bp}{\mathbb{P}}
\newcommand{\ba}{\mathbb{A}}
\newcommand{\spin}{{\rm Spin}\,}
\newcommand{\ds}{\displaystyle}
\newcommand{\tor}{{\rm Tor}\,}
\newcommand{\bff}{{\bf F}}
\newcommand{\bs}{\mathbb{S}}
\def\ns{\mathop{\lr}}
\def\nssup{\mathop{\lr\,sup}}
\def\nsinf{\mathop{\lr\,inf}}
\renewcommand{\phi}{\varphi}
\newcommand{\tT}{{\widetilde{T}}}
\newcommand{\tG}{{\widetilde{G}}}
\newcommand{\tB}{{\widetilde{B}}}
\newcommand{\tC}{{\widetilde{C}}}
\newcommand{\tW}{{\widetilde{W}}}
\newcommand{\tphi}{{\widetilde{\Phi}}}

\title{$K$-theory of Springer Varieties}

\author[P. Sankaran]{Parameswaran Sankaran}
\address{Chennai Mathematical Institute
H1, SIPCOT IT Park, Siruseri
Kelambakkam 603103
India}
\email{sankaran@cmi.ac.in}

\author[V. Uma]{Vikraman Uma}
\address{Department of Mathematics, Indian Institute of Technology, Madras, Chennai 600036, India
}
\email{vuma@iitm.ac.in}

\makeatletter
\@namedef{subjclassname@2020}{\textup{2020} Mathematics Subject Classification}
\makeatother

\begin{abstract}
  The aim of this paper is to describe the topological $K$-ring, in
  terms of generators and relations, of a Springer variety
  $\mathcal{F}_{\lambda}$ of type $A$ associated to a nilpotent
  operator having Jordan canonical form whose block sizes form a
  weakly decreasing sequence $\lambda=(\lambda_1,\ldots, \lambda_l)$.
  Our description parallels the description of the integral cohomology
  ring of $\mathcal{F}_{\lambda}$ due to Tanisaki and also the
  equivariant analogue due to Abe and Horiguchi.
  \end{abstract}

\

\maketitle
\section{Introduction}\label{flag}

Fix a positive integer $n$ and consider the complete flag
variety $\mathcal F(\mathbb C^n)$ (or more briefly $\mathcal F$)
defined as
\[\mathcal{F}(\mathbb{C}^n):=\{\underline{ V}:=(0=V_0\subset V_1\subset \cdots\subset V_n=\mathbb C^n)\mid
  \dim V_i=i ~~\mbox{for all}~~ i\}.\]

 Let $N:  \mathbb{C}^n\lra \mathbb{C}^n$ 
denote a 
nilpotent linear transformation of $\mathbb C^n$.   
The
Springer variety of type $A$ associated to $N$ denoted by
$\mathcal{F}_{N}$ is the closed subvariety of $\mathcal{F}$
defined as 
\[\{\underline{V}\in \mathcal{F}~~\mid~~NV_i\subset V_{i-1}~~\mbox{for
    all}~~ 1\leq i\leq n\}.\] The Springer variety $\mathcal F_N$ is
seen to be the subvariety of $\mathcal F$ fixed by the action of the
infinite cyclic group generated by 
the unipotent element $U=I_n+N\in SL(n,\mathbb C)$.   Moreover, denoting by
$\lambda=(\lambda_1,\lambda_2,\ldots,\lambda_l)$ the partition of $n$ 
where the $\lambda_j$ are the sizes of the diagonal blocks of the
Jordan canonical form of $N$, the variety $\mathcal F_N$
depends, up to isomorphism, only on the partition $\lambda$.  This is so, because two
different choices of nilpotent transformations corresponding to the
same partition $\lambda$ are conjugates in $GL(n,\mathbb C)$.  For
this reason, we assume that $N$ itself is in the Jordan canonical
form:
$N=J_{\lambda}:=\textrm{diag}(J_{\lambda_1},\ldots,
  J_{\lambda_l})$ with $\lambda_1\ge \cdots\ge \lambda_l$ and denote
  the Springer variety $\mathcal F_N$ by $\mathcal F_\lambda$.  (Here
  $J_p=(a_{i,j}) \in M_p(\mathbb C)$ is the matrix where
  $a_{i,i+1}=1, 1\le i<p,$ and all other entries are zero.)  If
  $\lambda=(1,\ldots,1)$, then $N=0$ and we have
  $\mathcal F_{\lambda}=\mathcal F(\mathbb{C}^n)=\mathcal F$.  At the other
  extreme, when $\lambda=(n)$, $N$ is a regular nilpotent element and
  we see that $\mathcal F_{(n)}$ is the one-point variety consisting only
  of the standard flag
  $0=E_0\subset E_1\subset \cdots\subset E_n=\mathbb C^n$ where $E_j$
  is spanned by the standard basis vectors $e_1,\ldots, e_j$ for
  $ 1\le j\le n$.

  Note that $\mathcal F_\lambda$ is stable by the action of the
  algebraic torus $T^l_\mathbb C\cong (\mathbb C^*)^l$ contained in
  $GL(n,\mathbb C)$ consisting of all diagonal matrices which commute
  with $N$.  We shall denote by $T^l=(\mathbb S^1)^l$ the compact
  torus contained in $T^l_\mathbb C$.  Denoting the diagonal subgroup
  of $GL(n,\mathbb C)$ by $T^n_\mathbb C$, we have
  $(t_1,\ldots, t_n)\in T^n_\mathbb C$ belongs to $T^l_\mathbb C$ if
  and only if $t_{a_j+i}=t_{a_j+1}$ for $1\le i\le \lambda_{j+1}$
  where $a_j:=\lambda_1+\cdots+\lambda_j$, $1\le j\le l-1$ and
  $a_0=0$.

The variety $\mathcal F_\lambda$ was first studied by Springer (see
\cite{spr1}, \cite{spr2} and also \cite{hs}). In particular, Springer
showed that there is a natural action of the symmetric group
$\mathfrak{S}_n$ on the rational cohomology
$H^*(\mathcal{F}_{\lambda};\mathbb{Q})$ which is compatible with the
standard action of $\mathfrak{S}_n$ on $H^*(\mathcal{F},\mathbb{Q})$.
Moreover, the restriction homomorphism
$H^*(\mathcal{F};\mathbb{Z})\lra H^*(\mathcal{F}_{\lambda};\mathbb{Z})
$ induced by the inclusion $\mathcal{F}_{\lambda}\hra \mathcal{F}$, is
surjective (see \cite{hs}).  The variety $\mathcal F_\lambda$ is not
irreducible in general, but it is equidimensional.  The irreducible
components of $\mathcal F_\lambda$ are naturally labelled by the set
of standard tableaux of shape $\lambda$.  See \cite{spal}.  Under the
$\mathfrak{S}_n$-action, the $\mathfrak{S}_n$-module,
$H^{*} (\mathcal{F}_{\lambda};\mathbb{Q}))$ is isomorphic to the
representation $M_\lambda$ of $\mathfrak{S}_n$ induced from the
identity representation of the subgroup
${\mathfrak{S}}_{\lambda}={\mathfrak{S}}_{\lambda_1}\times
{\mathfrak{S}}_{\lambda_2}\times \cdots \times
{\mathfrak{S}}_{\lambda_l}\subset {\mathfrak{S}}_n$.  See \cite{mac}
and \cite[p. 204]{dp}.

De Concini and Procesi \cite{dp} gave a description of
$H^*(\mathcal{F}_{\lambda};\mathbb{C})$ as the coordinate ring of an
(unreduced) variety over $\mathbb C$ which we now describe.  Let
$\lambda^\vee$ denote the partition dual to $\lambda$.  Consider the
coordinate ring
$\mathbb C[\mathfrak t_\mathbb C\cap \bar{O}_{\lambda^\vee}]$ of the
(non-reduced) scheme
$\mathfrak t_\mathbb C\cap \bar{O}_{\lambda^\vee}$ (scheme theoretic
intersection), where
$\mathfrak t_\mathbb C=Lie(T_\mathbb C^n)\subset
\mathfrak{gl}(n,\mathbb C)=M_n(\mathbb C)$ and
$ \bar O_{\lambda^\vee}\subset M_n(\mathbb C)$ denotes the closure of
the orbit of $J_{\lambda^\vee}$ under the adjoint action of
$GL(n,\mathbb C)$.  De Concini and Procesi showed that
$H^*(\mathcal F_\lambda;\mathbb C)$ is isomorphic to the algebra
$\mathbb C[\mathfrak t_\mathbb C\cap \bar O_{\lambda^\vee}]$.

Tanisaki \cite{t} described $H^*(\mathcal F_\lambda;\mathbb C)$ as a quotient of a
polynomial ring over $\mathbb{C}$ by an ideal, which has come to be known as
the Tanisaki ideal.   Tanisaki's description in fact yields the integral cohomology 
ring of $\mathcal F_\lambda$. 
Recently, the $T^l$-equivariant cohomology algebra
$H^*_{T^l}(\mathcal F_\lambda;\mathbb Z)$ has been described by H. Abe
and T. Horiguchi.  It turns out that
$H^*_{T^l}(\mathcal F_\lambda;\mathbb Z)$ is the quotient of a
polynomial algebra over $H^*_{T^l}(pt;\mathbb Z)=H^*(BT^l;\mathbb Z)$
modulo an ideal, which is a natural generalization of the Tanisaki
ideal.
This presentation
recovers the presentation for the ordinary integral cohomology ring
via the forgetful map $H^*_{T^l}(\mathcal{F}_{\lambda};\mathbb{Z})\lra
H^*(\mathcal{F}_{\lambda};\mathbb{Z})$.

Our aim in this paper is to describe the topological $K$-ring of the
variety $\mathcal{F}_{\lambda}$, in terms of generators and relations,
using Tanisaki's description of the integral cohomology ring.  The
generators of the Tanisaki ideal admit topological interpretation in
terms of the Chern classes of certain naturally defined line bundles
over $\mathcal F_\lambda$.  We interpret these relations as a {\it
  consequence} of a relation in the topological $K$-ring
$K(\mathcal F_\lambda)$ among the line bundles over
$\mathcal F_\lambda$.

 Before stating the main result of this paper, we need the
following notations.

A non-increasing sequence
$\lambda=(\lambda_1,\ldots, \lambda_l)$ of positive integers where
$\displaystyle\sum_{1\le j\le l} \lambda_j=n$, will be identified with
the partition $(\lambda_1,\ldots, \lambda_n)$ where $\lambda_j=0$ for
$j>l$.  Let $1\le s\le n$ and denote by $\mathbf i=(i_1,\ldots, i_s)$
a strictly increasing sequence of integers between $1$ and $n$.  The
set of all such sequences will be denoted $W_{n,s}$.  Let
$p_{\lambda}(s):=\lambda_{n-s+1}+\cdots+\lambda_n, 1\le s\le n$. 

We denote by $\mathcal L_i$ the canonical line bundle over $\mathcal F(\mathbb C^n)$ whose fibre over 
a flag $\underline{V}$ is the vector space $V_i/V_{i-1}, 1\le i\le n$.  Let $L_i=\mathcal L_i|_{\mathcal F_\lambda}$.  

We now state the main theorem.  A more precise formulation as a quotient of a polynomial ring by 
the $K$-theoretic Tanisaki ideal will be given in \S4. 

\bth\label{maintheorem} The ring $K^0(\mathcal F_\lambda)$ is
generated by the classes $[L_i]$ for $1\le i\le n$, subject only to
the following
(generating) relations:\\
\[\gamma^{d}([L_{i_1}\oplus \cdots\oplus L_{i_s}]-s)=0\]
for $d\ge s+1-p_{\lambda^\vee}(s),~~\mathbf i=(i_1,\ldots, i_s)\in W_{n,s}, ~1\le s\le n$.  Moreover,
$K^1(\mathcal F_\lambda)=0$.  \eeth

For the definition of $\gamma$-operations in $K$-theory, see \cite[Chapter 12]{husemoller}.

The structure of $K(\mathcal F_\lambda)$ as an abelian group is 
easily obtained from the fact that $\mathcal F_\lambda$ admits an
{\it algebraic} cellular decomposition.  The existence of 
such a decomposition was established by Spaltenstein \cite{spal}, the total number of such cells
being equal to the number of $T^l$-fixed points in
$\mathcal{F}_{\lambda}$ which in turn equals ${n\choose \lambda}$.
Although these algebraic cells are topologically $\mathbb C^r$, in general they
do not yield the structure of a CW complex.  See \cite{tym} for an
explicit example.  However, the algebraic cell-structure allows one to
compute the integral cohomology groups.  In particular,
$\mbox{rank} (H^*(\mathcal{F}_{\lambda};\mathbb{Z}))={n\choose
  \lambda}$ and $H^{k}(\mathcal{F}_{\lambda};\mathbb{Z})=0$ for $k$
odd.
 In particular, $K^1(\mathcal F_\lambda)=0$
and $K^0(\mathcal F_\lambda)$ is a free abelian group of rank
${n\choose \lambda}$.  
The fact that $H^*(\mathcal F_\lambda;\mathbb Z)$ is generated in degree $2$
allows us to apply \cite[Lemma 4.1]{su1} to
$\mathcal{F}_{\lambda}$ to 
conclude that $K^0(\mathcal F_\lambda)$ is generated as a {\it ring} by the classes of line bundles. 
  
Using a geometrical argument, we shall show that the ideal of relations 
among the $[L_j]$ contains the elements described in the theorem.  To show that these are all the generating 
relations we use Tanisaki's description of the {\it ring} $H^*(\mathcal F_\lambda;\mathbb Z)$ and a purely algebraic argument exploiting the fact that the $K$-theoretic Tanisaki ideal admits a certain natural filtration.  

Similar approaches were applied in our  
papers \cite{su1}, \cite{su2} and
\cite{s}, to obtain a presentation of the $K$-ring of a
smooth projective toric variety, a quasitoric manifold, and a class of torus manifolds 
which include the class of smooth complete toric varieties.

\section{Cohomology of Springer varieties}

Let $\mathcal V_j$ be the subbundle of the trivial vector bundle
$\mathcal{F}\times \mathbb C^n$ whose fibre over the flag
$\underline{V}=(V_i)\in \mathcal{F}$ is just $V_j$.  Denote by
$\mathcal L_i$ the line bundle
$\mathcal V_i/\mathcal V_{i-1}, 1\le i\le n,$ on $\mathcal{F}$.  We
denote the first Chern class of the line bundle $\mathcal L_i$ by
$x_i\in H^2(\mathcal{F};\mathbb{Z}), 1\le i\le n$.  One has an exact
sequence of algebraic vector bundles
$0\to \mathcal V_{s-1}\hookrightarrow \mathcal V_s\to \mathcal L_s\to
0$, which leads to an isomorphism of {\it complex} vector bundles for
$1\le s\le n$:
\begin{equation} \label{vsassumoflinebundles}
\mathcal L_1\oplus \cdots\oplus \mathcal L_s\cong \mathcal V_s.
\end{equation} 
Since, $\mathcal V_n=n\epsilon$, the trivial vector bundle of rank $n$, we have 
\begin{equation} \label{sumoftautlinebundles}
\mathcal L_1\oplus \cdots\oplus \mathcal L_n\cong n\epsilon.
\end{equation}
It follows that the Chern polynomial of
$\displaystyle\bigoplus_{1\le i\le n} \mathcal L_j$ is trivial. That
is, $\displaystyle\prod_{1\le i\le n}(1+x_it)=1$ and so the $j$th
elementary symmetric polynomial $e_j(x):=e_j(x_1,\ldots, x_n)$
vanishes for $1\le i\le n$.  Borel has shown that
$H^*(\mathcal F;\mathbb Z)$ is generated by the
$c_1(\mathcal L_i)=x_i, 1\le i\le n$ and that the only generating
relations among the Chern classes of $\mathcal L_i$ are given by
$e_j(x)=0, 1\le j\le n$.  (See \cite{borel}.) 
Thus we have a presentation
\begin{equation}
H^*(\mathcal F;\mathbb Z)\cong \frac{\mathbb Z[x_1,\ldots, x_n]}{\langle e_j(x);1\le j\le n\rangle}
.\end{equation}

For an arbitrary partition
$\lambda=\lambda_1\ge \cdots\ge \lambda_l$ of $n$, the Springer
variety $\mathcal F_\lambda$ is naturally imbedded in $\mathcal F$.
The cohomology ring of the Springer variety $\mathcal F_\lambda$ has
been described by Tanisaki \cite{t} in terms of generators and
relations, in a way that generalizes the above description of
$H^*(\mathcal F;\mathbb Z)$.  It turns out that, although Tanisaki
considered cohomology with complex coefficients, his description, 
recalled below, is valid even when the coefficient ring is the
integers (see \cite{dp}).

We need the following notation.   

\bdefe\label{plambda}
 Let $[n]:=\{1,2,\ldots, n\}$.  We  define the function $p_\lambda:[n]\to [n]\cup\{0\}$ associated to a partition $\lambda$ of $n$ as follows:
\begin{equation}
p_\lambda(s)=\lambda_{n-s+1}+\cdots+\lambda_n,~1\le s\le n.
\end{equation} 
\edefe Thus $p_\lambda$ is a monotonically increasing function of
$s$. The function $p_{\lambda^\vee}$ associated to the dual partition
$\lambda^\vee$ is more relevant for us.  Recall that the dual
partition $\lambda^\vee $ is defined as
$\lambda^\vee=(\eta_1,\ldots, \eta_n)$ where
$\eta_j=\# \{i\mid \lambda_i\ge j\}$.  Writing $\lambda$ as
$1^{a_1}\cdot 2^{a_2}\cdots n^{a_n}$, where $a_j$ is the number of times
$j$ occurs in $\lambda$, we have $\eta_j=a_j+\cdots+a_n$ for all
$j\ge 1$.

We pause for an example.

\beg Let $n=20$, and,
$\lambda=(5,4,4,2,2,2,1)=1^{1}\cdot 2^{3}\cdot 3^{0}\cdot 4^{2}\cdot
5^{1}$.  Then $\lambda^\vee=(7,6,3,3,1) $ and $p_{\lambda^\vee}(s)=0$
for $1\le s\le 15$,
$p_{\lambda^\vee}(16)=1,~ p_{\lambda^\vee}(17)=4,~
p_{\lambda^\vee}(18)=7, ~p_{\lambda^\vee}(19)=13,~
p_{\lambda^\vee}(20)=20$.  \eeg

\bdefe\label{tanisakiideal} Let $S=\mathbb Z[y_1,\ldots,y_n]$ be the polynomial ring in $n$-indeterminates $y_1,\ldots, y_n$.  The {\em Tanisaki ideal} is the ideal 
$I_\lambda\subset S$ generated by the following elements:
\[e_d(y_{i_1},\ldots, y_{i_s}), ~\mbox{for}~ d\ge
  s+1-p_{\lambda^\vee}(s)\] where
$1\le i_1<\cdots<i_s\le n, ~1\le s\le n.$
 \edefe

Recall that $L_j$ is the line bundle over $\mathcal F_\lambda$
obtained as the restriction of $\mathcal L_j$ on $\mathcal F$.

\bth  \label{tanisaki}
We keep the above notations. 
Let $\lambda=\lambda_1\ge \cdots\ge \lambda_l$ be a partition of $n$.   Then one has an isomorphism of rings 
\[H^*(\mathcal F_\lambda;\mathbb Z)\cong  {S}/{I_\lambda }\]
where $c_1(L_j)$  corresponds to $y_j+I_\lambda, 1\le j\le n$.  Moreover,  the inclusion $\iota_\lambda: \mathcal F_\lambda
\to \mathcal F$ induces a surjection $\iota_\lambda^*: H^*(\mathcal F;\mathbb Z)\to H^*(\mathcal F_\lambda;\mathbb Z)$.\hfill $\Box$
\eeth
We shall abuse notation and denote by $y_j\in H^*(\mathcal F_\lambda;\mathbb Z)$ the image of $y_j+I_\lambda$ under the isomorphism of the above theorem. 

The rank of $H^*(\mathcal F_\lambda;\mathbb Z)$ equals
\[{n\choose \lambda}=\frac{n!}{(\lambda_1!\cdots \lambda_l!)}\] and
$\displaystyle\dim_\mathbb C\mathcal F_\lambda=\sum_{1\le j\le
  n}\eta_j.(\eta_j-1)/2$ where
$\lambda^\vee=\eta_1\ge\cdots\ge \eta_n\ge 0$ is the partition of $n$
that is dual to $\lambda$.

\section{Line bundles over $\mathcal F_\lambda$}

We begin by recalling the description of $K(\mathcal F)$ in terms of
generators and relations.  It turns out that $K(\mathcal F)$ is
generated by the classes of the line bundles
$[\mathcal L_j], 1\le j\le n$.  (This can be seen using the
Atiyah-Hirzebruch spectral sequence and the Chern character map
$\textrm{ch}:K(\mathcal F) \to H^*(\mathcal F;\mathbb Q)$. (See
\cite{athir}).) In view of (\ref{sumoftautlinebundles}), we have an
isomorphism
\begin{equation}\label{kringflag}
e_k(\mathcal L_1, \ldots,\mathcal L_n)\cong {n\choose k}\epsilon
\end{equation} for all $k\ge 1$.  
Indeed, we have an isomorphism 
\begin{equation}
\phi: \frac{\mathbb Z[u_j;1\le j\le n]}{\mathcal I} \to K(\mathcal F)
\end{equation}
where $\phi(u_i )=[\mathcal L_i], 1\le i\le n,$ and the ideal $\mathcal I$ is generated by the elements 
$e_k(u_1,\ldots, u_n)-{n\choose k}, 1\le k\le n$.   (See \cite[\S3, Chapter IV]{k}.)

Recall that the Atiyah-Hirzebruch spectral sequence of a space $X$
admitting the structure of a finite CW complex has $E_2$-page defined
as $E_2^{p,q}(X)=H^p(X;K^q(pt))$ and differential $d_r$ has bidegree
$(r,1-r)$.  The spectral sequence converges to $K^*(X)$, that is,
there exists a decreasing filtration $\{K_p^{p+q}(X)\}$ such that
$E_\infty^{p,q}\cong G_pK^{p+q}(X)=K_{p}^{p+q}(X)/K_{p+1}^{p+q}(X)$
(see \cite[p. 17]{athir}).

Recall that $K^{q}(pt)\simeq \mathbb{Z}$ for $q$ even and
$K^{q}(pt)=0$ for $q$ odd (see \cite[p. 10]{athir}).  Suppose that
$H^p(X;\mathbb Z)$ vanishes if $p$ is odd and that it is a free
abelian group when $p$ is even. Then $E_2^{p,q}=0$ unless both $p,q$
are even and it follows that $d_r=0$ for all $r$. Thus the spectral
sequence collapses and we have $E_2^{p,q}=E^{p,q}_\infty.$
Consequently, $K^1(X)=0$ and $K^0(X)\cong H^*(X;\mathbb Z)$ is a free
abelian group.

Suppose that $Y$ is another such space and that $f:X\to Y$ induces a
surjective homomorphism $f^*:H^*(Y;\mathbb Z)\to H^*(X;\mathbb Z)$.
The naturality of the spectral sequence yields a morphism
$\{E_r^{p,q}(Y),d_r\}\to \{E_r^{p,q}(Y),d_r\}$.  Since
$E_2^{p,q}(Y)=H^p(Y;K^q(pt)) \stackrel{f^*}{\longrightarrow }
H^p(X;K^q(pt))=E_2^{p,q}(X)$ is surjective and since the differentials
$d_r$ vanish for both the sequences for $r\ge 2$, we see that
$f^*:E_\infty^{p,q}(Y)=E_2^{p,q}(Y)\to E_2^{p,q}(X
)=E_\infty^{p,q}(X)$ is surjective.  Moreover, both
$K^*(X)=K^0(X),K^*(Y)=K^0(Y)$ are free abelian groups.  It follows
that $f^!:K^*(Y)\to K^*(X)$ is surjective.

Setting $X=\mathcal F_\lambda$, we obtain the following isomorphism of
abelian groups
\begin{equation}
K^*(\mathcal F_\lambda)=K^0(\mathcal F_\lambda)\cong H^*(\mathcal F_\lambda;\mathbb Z)\cong \mathbb Z^{{n\choose \lambda}}.
\end{equation}

We next take $Y=\mathcal{F}$ and $f$ to be the inclusion
$\iota_\lambda:\mathcal F_\lambda\hookrightarrow \mathcal F$. 

Let $y_i$ denote the image
$\iota^*_{\lambda}(x_i)\in H^2(\mathcal F_\lambda;\mathbb Z)$ and
$L_i$ denote the restriction $\mathcal L_i|_{\mathcal F_\lambda}$ for
$1\leq i\leq n$.  Thus $c_1(L_i)=y_i$.  The surjectivity of
$\iota_\lambda^*:H^*(\mathcal F;\mathbb Z)\to H^*(\mathcal
F_\lambda;\mathbb Z)$ implies that $H^*(\mathcal F_\lambda;\mathbb Z)$
is generated by $c_1(L_j)=y_j ; 1\le j\le n$.

By the above argument we conclude that the pull back
map \[\iota_\lambda^!:K(\mathcal F)\to K(\mathcal F_\lambda)\] induced
by the inclusion $\iota_\lambda:\mathcal F_\lambda \to \mathcal F$ is
surjective.  Furthermore, since $K(\mathcal F)$ is generated, as a
ring by the isomorphism classes of line bundles
$[\mathcal L_j], 1\le j\le n$, it follows that $K(\mathcal F_\lambda)$
is generated by $[L_j],1\le j\le n$.

Summarising the above discussion, we have the following proposition.

\bpropo \label{linebundlesgenerate} The group $K(\mathcal F_\lambda)$
is free abelian of rank ${n\choose \lambda}$. Also,
$\iota_\lambda^!:K(\mathcal F)\to K(\mathcal F_\lambda)$ is a
surjection.  In particular, $K(\mathcal F_\lambda)$ is generated as a
ring by the classes of line bundles
$[L_j]\in K(\mathcal F_\lambda),1\le j\le n$.  \hfill $\Box$ \epropo

\subsection{The action of the symmetric group on $K(\mathcal F_\lambda)$}
The symmetric group $\mathfrak{S}_n$ acts linearly on $\mathbb C^n$ by
permuting the standard basis vectors $e_1,\ldots, e_n$.  This action
induces an action of $\mathfrak{S}_n$ on the flag variety $\mathcal F$
and hence on $K(\mathcal F) $ as well as on the singular cohomology
algebra $H^*(\mathcal F; \mathbb Z)$.  Explicitly, the
$\mathfrak{S}_n$-action on $K(\mathcal F)$ and on
$H^*(\mathcal F;\mathbb Z)$ are given by permutation of the classes of
line bundles $[\mathcal L_j]\in K(\mathcal F)$ and the Chern classes
$x_j:=c_1(\mathcal L_j)\in H^*(\mathcal F;\mathbb Z), 1\le j\le n$,
respectively.  It is readily seen that the Chern character map
\[\textrm{ch}_\mathcal F:K(\mathcal F)\to H^*(\mathcal F;\mathbb Q)\] which sends $[\mathcal L_k]$ to 
$\displaystyle \sum_{r\ge 0}\frac{x_k^r}{r!}$ is $\mathfrak{S}_n$-equivariant.

Springer \cite{spr2} showed that the symmetric group $\mathfrak{S}_n$ acts on $H^*(\mathcal F_\lambda;\mathbb Q)$ and that $\iota^*_\lambda:H^*(\mathcal F;\mathbb Q)\to H^*(\mathcal F_\lambda;\mathbb Q)$ is $\mathfrak{S}_n$-equivariant and surjective.  Hotta and Springer \cite{hs} showed that in fact $\iota^*_\lambda:H^*(\mathcal F;\mathbb Z)\to H^*(\mathcal F_\lambda;\mathbb Z)$ 
is $\mathfrak{S}_n$-equivariant and surjective (also see \cite[p. 213]{dp}).

\bpropo \label{snaction} The ring $K(\mathcal F_\lambda)$ admits an
action of the symmetric group $\mathfrak{S}_n$ with respect to which the Chern
character map
\[\textrm{ch}_{\mathcal F_\lambda}:K(\mathcal F_\lambda)\to
  H^*(\mathcal F_\lambda;\mathbb Q)\] as well as the pull back map
$\iota^!_\lambda: K(\mathcal F)\to K(\mathcal F_\lambda)$ are both
$\mathfrak{S}_n$-equivariant.  \epropo
\begin{proof}
By the naturality of Chern character we have the following commuting diagram:
\begin{equation}\label{cherncomm}
\begin{array}{ccc}
K(\mathcal F)&\stackrel{\iota_\lambda^!}{\longrightarrow}& K(\mathcal F_\lambda)\\
\textrm{ch}_{\mathcal F}\downarrow &&\downarrow \textrm{ch}_{\mathcal F_\lambda}\\
H^*(\mathcal F;\mathbb Q)&\stackrel{\iota^*_\lambda}{\longrightarrow}& H^*(\mathcal F_\lambda;\mathbb Q).\\
\end{array}
\end{equation}
The vertical arrows are monomorphisms since the $K$-groups
$K(\mathcal F), K(\mathcal F_\lambda)$ are free abelian groups and
$\textrm{ch}\otimes \mathbb Q$ is an isomorphism.

Now, since $\textrm{ch}_{\mathcal{F}_{\lambda}}$ is injective,
$\iota_{\lambda}^!$ is surjective and $\iota_\lambda^*$ and
$\textrm{ch}_\mathcal F$ are $\mathfrak{S}_n$-equivariant, the
commutative diagram \eqref{cherncomm} implies that the image
$\textrm{ch}_{\mathcal{F}_{\lambda}}(K(\mathcal F_\lambda))$ in
$H^*(\mathcal{F}_{\lambda};\mathbb{Q})$ is stable under the
$\mathfrak{S}_n$-action. This further implies that
$K(\mathcal F_\lambda)$ admits an action of $\mathfrak{S}_n$ such that
$\textrm{ch}_{\mathcal F_\lambda}$ and $\iota_\lambda^!$ are
equivariant.  This can be seen more explicitly as follows.

Let $x\in K(\mathcal{F}_{\lambda})$. Since $\iota_{\lambda}^!$ is
surjective, $x=\iota_{\lambda}^!(y)$ for some $y\in
K(\mathcal{F})$. We define the action of $\mathfrak{S}_n$ on (the
right of) $K(\mathcal{F}_{\lambda})$ as follows:
$ x\cdot\sigma:=\iota_{\lambda}^!(y\cdot\sigma)$. Indeed from
\eqref{cherncomm} we have
\[\textrm{ch}_{\mathcal{F}_{\lambda}}\circ \iota_{\lambda}^!(y\cdot\sigma)=\iota_{\lambda}^*\circ \textrm{ch}_{\mathcal{F}} (y \cdot\sigma).\]
Moreover, since $\iota_{\lambda}^*$ and $\textrm{ch}_{\mathcal{F}}$
are $\mathfrak{S}_n$-equivariant this implies
\[\textrm{ch}_{\mathcal{F}_{\lambda}}\circ \iota_{\lambda}^!(y \cdot\sigma)=
  (\iota_{\lambda}^*\circ \textrm{ch}_{\mathcal{F}}(y))\cdot\sigma.\] Again by
the commutativity of \eqref{cherncomm} we get
\[\textrm{ch}_{\mathcal{F}_{\lambda}}\circ \iota_{\lambda}^! (y\cdot\sigma)=
  (\textrm{ch}_{\mathcal{F}_{\lambda}}\circ \iota_{\lambda}^!(y))
  \cdot\sigma= \textrm{ch}_{\mathcal{F}_{\lambda}}(x)\cdot\sigma.\]
Thus by definition
$\textrm{ch}_{\mathcal{F}_{\lambda}}(x\cdot\sigma)=
\textrm{ch}_{\mathcal{F}_{\lambda}}(x)\cdot\sigma$. Since
$\textrm{ch}_{\mathcal{F}_{\lambda}}$ is injective and the
$\mathfrak{S}_n$ action on $H^*(\mathcal{F}_{\lambda};\mathbb{Q})$ is
well defined, it follows that the $\mathfrak{S}_n$ action on
$K(\mathcal{F}_{\lambda})$ is well defined. It also follows that
$\textrm{ch}_{\mathcal{F}_{\lambda}}$ and $\iota^{!}_{\lambda}$ are
$\mathfrak{S}_n$-equivariant.
\end{proof}

Since the action of $\mathfrak S_n$ on (the right of) $K(\mathcal F)$
is obtained as $[\mathcal L_j]\cdot \sigma = [\mathcal L_{\sigma(j)}]$
for $1\le j\le n,$ and $\sigma \in \mathfrak S_n$, by Proposition
\ref{snaction} and Proposition \ref{linebundlesgenerate}, it follows
that the $\mathfrak S_n$ action on $K(\mathcal F_\lambda)$ is given by
$[L_j]\cdot\sigma=[L_{\sigma(j)}]$ for all
$1\le j\le n, ~~\sigma\in \mathfrak S_n$.

Recall the function $p_\lambda$ defined in \S2. 
The numbers $p_{\lambda^\vee}(s)$ are related to the nilpotent
transformation $N=J_{\lambda}$ as follows. 

\blem\label{rank} \cite[Proposition 3]{t}.
With the above notations, \\(i) $p_{\lambda^\vee}(s)=\mathrm{rank}(J_\lambda^{n-s}), ~1\le s\le n$.\\
(ii) Let $\underline U=U_1\subset U_2\subset \cdots\subset U_n=\mathbb C^n$ be a flag that refines the partial flag 
$0=\im (J_\lambda^{\lambda_1})\subset \im(J_\lambda^{\lambda_1-1})\subset\cdots \subset \im (J_\lambda^2)\subset \im(J_\lambda)\subset \mathbb C^n$.
Then, for any $\underline V\in \mathcal F_\lambda$  and any $s\ge 1$, we have  
$U_q\subset V_s$  where  $q=p_{\lambda^\vee}(s)$.  \hfill $\Box$
\elem

\subsection{Sectioning canonical bundles over $\mathcal F_\lambda$}
For $1\le s\le n$ we let \begin{equation}\label{W_{n,s}} W_{n,s}:=\{\mathbf
  i=(i_1,\ldots, i_s) ~~\mbox{where}~~ 1\le i_1<\cdots<i_s\le
  n\}. \end{equation}

\bpropo \label{relationsinflambda} Let $1\le s\le n$ and let
$\mathbf {i}\in W_{n,s}.$ Then
\[L_{i_1}\oplus \cdots \oplus L_{i_s}\cong \xi\oplus q\epsilon\] for some complex vector bundle $\xi=\xi(\mathbf{i})$ over $\mathcal F_\lambda$ where 
$q:=p_{\lambda^\vee}(s)$.  
\epropo
\begin{proof}  Fix $s\le n$.   Since the action of $\mathfrak{S}_n$ on
  $K(\mathcal F_\lambda)$ permutes the $L_j, ~1\le j\le n$, we need
  only consider the case where $i_j=j, ~~\forall ~~1\leq j\leq s$.

We replace $N$ by a conjugate $gN g^{-1}$ 
so that $\im (gNg^{-1})^{n-k}=U_{p_{\lambda^\vee}(k)}=\mathbb C^{p_{\lambda^\vee}(k)}$ for $k\ge 1$.  
We may then choose the refinement $\underline{U}\in \mathcal F$ to be the standard flag 
$0\subset \mathbb C\subset \cdots\subset\mathbb C^n$.   
Thus $\mathbb C^q\subset V_s$ for any $\underline{V}\in \mathcal F_{gNg^{-1}}$.  Let $\iota_g:\mathcal F\to \mathcal F$ be 
the translation by $g$: $\underline {V}\mapsto g\underline {V}=gV_0=0\subset gV_1\subset \cdots \subset gV_n=\mathbb C^n$.  
Since $GL(n,\mathbb C)$ is connected, the composition $\mathcal F_{N}\stackrel{\iota_\lambda}{\hookrightarrow} \mathcal F\stackrel{\iota_g}{\longrightarrow }\mathcal F$, denoted $\iota_{\lambda,g}$ is homotopic to $\iota_\lambda$ and maps $\mathcal F_{N}$ 
onto $\mathcal F_{gNg^{-1}}\subset \mathcal F$.  It follows that
$\iota_\lambda$ and $\iota_{\lambda_g}$ induce the same map in
$K$-theory and singular cohomology.  In particular,
$\iota_{\lambda,g}^*(\mathcal L_j)=L_j~\forall ~1\leq j\le n$.

Let $G_{n,s}=G_s(\mathbb C^n)$ denote the Grassmann variety of $s$-planes in $\mathbb C^n$.   One has 
a projection $\pi_s: \mathcal F_{N}\to G_{n,s}$ defined as $\underline{V}\mapsto V_s$.  

Let $Y_q\subset G_{n,s} $ denote the subvariety $\{U\in G_{n,s}\mid U\supset U_q\}, 1\le q<s$.  
Then $Y_q$ is isomorphic to a Grassmann variety $G_{n-q,s-q}$. A specific isomorphism  $Y_q\cong G_{s-q}(\mathbb C^n/U_q)$ is
obtained by sending $U\in Y_q$ to $U/U_q$.   The tautological complex vector bundle $\gamma_{n,s}$ is of rank $s$, whose fibre over $A\in G_{n,s}$ is the vector space $A$.  When restricted to $Y_q$, $\gamma_{n,s}$  has a trivial subbundle $q\epsilon $ of rank $q$.  Indeed we have a commuting diagram 
\[\begin{array}{ccc}Y_q\times U_q& \hookrightarrow&  E(\gamma_{n,s}|_{Y_q})\\
\downarrow && \downarrow \\
Y_q&\stackrel{id}{\longrightarrow}& Y_q\\
\end{array}\]
where the vertical arrows are bundle projections. 
Therefore
\begin{equation} \label{splitting}
\gamma_{n,s}|_{Y_q}\cong \omega\oplus q\varepsilon
\end{equation}
where $\omega$ is the complex vector bundle over $Y_q$ whose fibre over $A\in Y_q$ is the complex vector 
space $A':=A/\mathbb C^q.$

From Proposition \ref{rank}, the image of the composition
\[\mathcal{F}_{N}=\mathcal F_\lambda \stackrel{\iota_{\lambda,g}}{\longrightarrow}\mathcal F\stackrel{\pi_\lambda}{\longrightarrow} G_{n,s},\] denoted
$\pi_{\lambda,s}$, is contained in $Y_q$.  Therefore, we have a
commuting diagram
\begin{equation}
\begin{array}{rcc}
\mathcal F_\lambda &\stackrel{\iota_{\lambda,g}}{\longrightarrow} & \mathcal F\\
\pi_{\lambda,s}\downarrow &&~~\downarrow \pi_s\\
Y_q& \hookrightarrow & G_{n,s}.\\
\end{array}
\end{equation}

Now
$\pi_s^*(\gamma_{n,s})=\mathcal V_s=\mathcal L_1\oplus \cdots \oplus
\mathcal L_s$ by (\ref{vsassumoflinebundles}).  Therefore
$L_1\oplus \cdots\oplus L_s=\iota_{\lambda,g}^*(\mathcal
L_1\oplus\cdots\oplus \mathcal L_s) =\iota_{\lambda,g}^*\circ
\pi_s^*(\gamma_{n,s})=\pi_{\lambda,s}^*(\gamma_{n,s}|_{Y_q})=\pi_{\lambda,s}^*(\omega)\oplus
q\epsilon$, from (\ref{splitting}).
\end{proof}

{\it Remark:} Note here that $Y_q$ is nothing but the Schubert variety
\[X(\sigma)=\{U\in G_{n,s}\mid \dim(U\cap \mathbb C^{\sigma_i})\ge
  i,~~ 1\le i\le s\}\] where
\[\sigma_i=\begin{cases}i, & \textrm{if~} i\le q,\\
  n-s+i,&\textrm{if~} q<i\le s.
\end{cases}\]

\subsection{The $\gamma$-operations in $K$-theory.}   We recall here  the $\gamma$-operations in $K$-theory and their relation to the $\lambda$-operations.  We refer the reader to  \cite[\S3, Chapter 12]{husemoller} for further details.
 Let $X$ be a finite CW complex.   The $\gamma$-operations 
$\gamma^d:K(X)\to K(X), d\ge 0,$ are defined in terms of the exterior power operations $\lambda^j$ as follows: \\
\begin{equation}\label{gamma-defn}
\gamma_t(x)=\lambda_{t/(1-t)}(x)~~\forall x\in K(X),
\end{equation}
where
$\displaystyle\gamma_t(x)=\sum_{k\ge 0}
  \gamma^k(x)t^k$ and
$\displaystyle\lambda_t(x)=\sum_{k\ge 0} \lambda^k
  (x)t^k$ are regarded as elements of the formal power series ring
$K(X)[[t]]$ in the indeterminate $t$.  Since
$\gamma^0(x)=\lambda^0(x)=1$, we can express $\lambda_t(x)$ in terms
of $\gamma(x)$. Indeed we have $\lambda_t(x)=\gamma_{t/(1+t)}(x)$.  It
follows from the definition of $\gamma^d(x)$ that
\begin{equation}\label{gammalambda}
\displaystyle \gamma^d(x)=\sum_ {0\le k\le d}\lambda^k(x){k+d-1\choose k-1}=\lambda^d(x+d-1).
\end{equation}
When $x=[\xi]\in K(X)$ is the class of  a vector bundle $\xi$  of rank
$k$, we have $\lambda_t(x)$ is a polynomial of degree $k$ since the
exterior power
\begin{equation}\label{lambdaequality}\lambda^d([\xi])=0=\lambda^d([\xi]+d-1-k)\end{equation}
for $d\ge k+1$.

If $\xi$ is a line bundle, then $\lambda_t(x)=1+xt$.  In the case when
$\xi$ is a trivial bundle, we have $[\xi]=k\in K(X)$ and
$\lambda_t(k)=(1+t)^k$ and so $ \lambda_t(-k)=(1+t)^{-k}$. The last
equality in \eqref{lambdaequality} follows from the identity
$\lambda_t(x+y)=\lambda_t(x)\lambda_t(y)$.

These basic facts will be used below.

\section{Proof of Theorem \ref{maintheorem}}
We begin by establishing the following proposition, which is an 
an immediate corollary of Proposition  \ref{relationsinflambda}.

\bpropo \label{mainrelations}
For any $\mathbf i\in W_{n,s}$ and any $d\ge s+1-p_{\lambda^\vee}(s)$, the following relation holds in 
$K(\mathcal F_\lambda)$:
\begin{equation} \label{kringspringer}
\gamma^d([L_{i_1}]+ \cdots+ [L_{i_s}]-s)= \lambda^d([L_{i_1}]+\cdots+[L_{i_s}]+d-s-1)=0
\end{equation}
for all $d\ge s+1-p_{\lambda^\vee}(s).$  Equivalently,
\begin{equation}
  \displaystyle\lambda^d((\sum_{1\le j\le s}[
  L_{i_j}])-q)=0 ~~\forall ~~d\ge s+1-q
\end{equation}  where
$q=p_{\lambda^\vee}(s)$.
\epropo

\begin{proof} Set
  $\displaystyle\mathcal E=\bigoplus_{1\le j\le s} L_{i_j}$.  Then, by
  Proposition \ref{relationsinflambda},
  $\displaystyle\mathcal E\cong \xi+q\epsilon$ where rank of $\xi$
  equals $s-q$.  (Here $q=p_{\lambda^\vee}(s)$ as in the Propositon.)
  It follows that
  $\displaystyle \gamma^d([\mathcal E]-s)=\gamma^d([\xi]-s+q)=
  \lambda^{d}([\xi]+d-s+q-1)=0$ for all
  $d\ge \textrm{rank}(\xi)+1=s+1-q$.  The last equation is equivalent
  to the vanishing of $\lambda^d([\mathcal E]-q)$ for
  all $d\ge s+1-q$ (see \eqref{lambdaequality}).
\end{proof}

Let $\mathbf i\in W_{n,s}$.  We have 
\begin{equation}\label{genrel1}\displaystyle\lambda_t(\sum_{1\le j\le
    s}[L_{i_j}]-q)=\prod_{1\le j\le
    s}(1+[L_{i_j}]t)\cdot (1+t)^{-q}\end{equation}
is a polynomial of degree $s-q$ where $q=p_{\lambda^\vee}(s)$.
Comparing the coefficient of $t^{d}$ where $d\ge s+1-q$,  
we obtain the following equation in $K(\mathcal F_\lambda)$ from Proposition
\ref{mainrelations}:
\begin{equation} \label{generatingrelations}
h_d([L_{i_1}],\ldots, [L_{i_s}])
:=\sum_{0\le k\le d} (-1)^{d-k}e_k([L_{i_1}],\ldots, [L_{i_s}])\cdot {q+d-k-1\choose q-1}=0.
\end{equation} 

When $s=n$, we have $q=n$. Now, from Equation \ref{sumoftautlinebundles} we have 
$\displaystyle\sum_{1\le j\le n}[L_{i_j}]-n=0$. Using this  in  Equation (\ref{genrel1}), we
get $\displaystyle (1+t)^{n}=\prod_{1\le j\le n}(1+[L_{i_j}]t)$. Thus
in this case (\ref{kringspringer}) follows from (\ref{kringflag}).

Let  $R:= \mathbb Z[u_1,\ldots,u_n]$, the polynomial algebra in $n$ variables $u_1,\ldots, u_n$. Let $1\le s\le n$ and let $\mathbf i\in W_{n,s}.$ 
The ring $R$ is graded by setting $\deg(u_j)=1~~\forall~~ j$.  One has the augmentation $\theta: R\to \mathbb Z$ defined  
by $u_j\mapsto 1~~\forall~~ j\le n$.   
For $\mathbf i=(i_1,\ldots, i_s)\in W_{n,s},$  we denote by $u_{\mathbf i}$ the sequence $(u_{i_1},\ldots, u_{i_s})$.
We define the {\it K-theoretic Tanisaki ideal}
$\mathcal I_\lambda\subset R$ to be the ideal generated by the
elements $h_d(u_{\mathbf i})$ defined as
\[ h_d(u_{\mathbf i}):=\sum_{0\le k\le d} (-1)^{d-k}e_k(u_{\mathbf
    i})\cdot {q+d-k-1\choose q-1},\] where
$1\le i_1<\ldots\le i_s\le s, ~1\le s\le n$.  Since
\[\sum_{0\le k\le d}(-1)^{d-k}{s\choose k}{q+d-k-1\choose
  q-1}={s-q\choose d}=0\] as $d\ge s+1-q,$ we have
\[h_d(u_{\mathbf i})=\sum_{0\le k\le d} (-1)^{d-k}(e_k(u_{\mathbf
    i})-{s\choose k})\cdot {q+d-k-1\choose q-1}.\]

The highest degree term in $h_d(u_\mathbf i)$ equals
$e_d(u_\mathbf i)$.  We set
\[\bar h_d(u_\mathbf i):=e_d(u_\mathbf i)-{s\choose
      d}\] so that
$\theta (\bar{h}_d(u_\mathbf i))=0=\theta(h_d(u_\mathbf i))$.

We may restate Theorem \ref{maintheorem} in view of Proposition \ref{mainrelations} and Equations (\ref{genrel1}) and (\ref{generatingrelations}) as follows:\\

\bth \label{maintheorem}  
We keep the above notations.   Let
\[\psi_\lambda:R=\mathbb{Z}[u_1,\ldots,u_n]\to
  K(\mathcal F_\lambda)\] be the ring homomorphism defined by
$\psi_\lambda(u_j)=[L_j], ~~1\le j\le n$.  Then $\psi_{\lambda}$ is
surjective and $\ker(\psi_\lambda)=\mathcal I_\lambda$.
 \eeth
\begin{proof}
  In view of (\ref{generatingrelations}), it is clear that
  $\mathcal I_\lambda\subset \ker(\psi_\lambda)$.  Also, by
  Proposition \ref{linebundlesgenerate}, $\psi_{\lambda}$ is
  surjective.  Since
  $\textrm{rank}(K(\mathcal F_\lambda))=\textrm{rank}(H^*(\mathcal
  F_\lambda;\mathbb Z))={n\choose\lambda},$ it suffices to show that
  $R/\mathcal I_\lambda$ is a free abelian group of rank
  ${n\choose \lambda}$.  We will show that $R/\mathcal I_\lambda$ has
  a filtration such that the associated graded ring is isomorphic as
  an abelian group to
  $S/I_\lambda\cong H^*(\mathcal F_\lambda;\mathbb Z)$.

  Let $R^d\subset R$ denote the abelian group generated by elements of
  the form $P(u)-\theta(P(u))\in \ker(\theta)$, where $P(u)$ is a
  polynomial in $u_1,\ldots, u_n$ of degree at most $d$.  Note that
  $R^0=0$ and $R^d\subset R^{d+1}$ and
  $R^i\cdot R^j\subset R^{i+j} ~~\forall~~ i, j \ge 1$.  Clearly
  $h_d(u_\mathbf i)-\bar h_d(u_\mathbf i)\in R^{d-1}$.

  Let $k\ge 1.$ and let $\mathbf s=(s_1,\ldots,s_k)$ where
  $1\le s_1\le \cdots \le s_k\le n$.  We denote by $\mathcal S_k$ the
  set of all such $k$-tuples $\mathbf s$.  For
  $\mathbf s\in \mathcal S_k$ we define the following sets
\[\mathcal W_{\mathbf s}:= W_{n,s_1}\times \cdots\times
  W_{n,s_k} ~~~\mbox{and}\]
\[\mathcal D_{\mathbf s}:=\{\mathbf d=(d_1,\ldots, d_k)\mid
  s_j+1-p_{\lambda^\vee}(s_j)\le d_j\le s_j~~\forall~~ j\}.\] For a
given $\mathbf d\in \mathcal D_{\mathbf s}$ and an integer
$\displaystyle d\ge \sum_{j=1}^k d_j$, we set
\[\mathcal C_{d,\mathbf d}:=\{\mathbf c=(c_1,\ldots, c_k)\mid c_j\ge
  1, \sum_{1\le j\le k} c_j\cdot d_j=d\}.\] An element of $\mathcal W_{\mathbf s}$ will be denoted $\iota=( \mathbf i(1),\ldots,\mathbf i(k))$.  

Given
$\mathbf s\in \mathcal S_k, ~\iota\in \mathcal W_\mathbf s, ~{\mathbf
  d}\in \mathcal D_\mathbf s, ~\mathbf c\in
\mathcal C_{d, \mathbf d}$, we set
\[h_{\iota,\mathbf d, \mathbf c}(u):= h_{d_1}(u_{\mathbf
    i(1)})^{c_1}\cdots h_{d_k}(u_{\mathbf i(k)})^{c_k}\in \mathcal
  I_\lambda^d\] where
\[\mathcal I^d_\lambda=\mathcal I_\lambda\cap R^d.\] Let
$\mathcal H_{d,k}$ be the set of all elements
$h_{\iota, \mathbf d,\mathbf c}(u),$ as we vary
$ \mathbf s\in\mathcal S_k, \iota\in \mathcal W_\mathbf s, \mathbf
d\in \mathcal D_\mathbf s, \mathbf c\in \mathcal C_{d, \mathbf d}$.
Let \[\mathcal H_d=\bigcup_{k} \mathcal H_{d,k}.\]

Note that the leading term of $h_{\iota, \mathbf {d},\mathbf c}(u)$,
when expressed as a polynomial in $u_1,\ldots, u_n,$ equals that of
$e_{\iota, \mathbf d, \mathbf c}(u):=e_{d_1}(u_{\mathbf
  i(1)})^{c_1}\cdots e_{d_k}(u_{\mathbf i(k)})^{c_k}$.  Since
$\theta (e_{\iota,\mathbf d,\mathbf c}(u))\ne 0,$ the element
$e_{\iota,\mathbf d,\mathbf c}(u)$ is not in
$R^d$.  But evidently
$e_{\iota,\mathbf d,\mathbf c}(u)-\theta(e_{\iota,\mathbf d,\mathbf
  c}(u))$ belongs to $R^d$.  Setting $v_j=u_j-1, 1\le j\le n,$ we have
$e_{\iota, \mathbf d,\mathbf c}(v)\equiv e_{\iota,\mathbf d,\mathbf
  c}(u)-\theta( e_{\iota,\mathbf d,\mathbf c}(u))\equiv h_{\iota,
  \mathbf d, \mathbf c}(u) \mod R^{d-1}$.

Recall, from Definition \ref{tanisakiideal}, that $I_\lambda\subset S=\mathbb Z[y_1,\ldots, y_n]$ is the
Tanisaki ideal generated by the homogeneous polynomials
\[e_{d}(y_{\mathbf i}), ~\mathbf i\in \mathcal W_{n,s}, ~d\ge
  s+1-p_{\lambda^\vee}(s), ~1\le s\le n\] (see Definition
\ref{tanisakiideal}).  Let $S_k$ denote the homogeneous polynomials in
$S$ of degree $k$ where $\deg(y_j)=1$ for all $j\le n$ and
$\displaystyle I_{\lambda,k}:=I_\lambda\cap S_k$.  Writing
$\displaystyle S^k=\bigoplus_{1\le j\le k}S_j$ and
$\displaystyle I^k_\lambda=\bigoplus_{1\le j\le k}I_{\lambda,j}$, we
see that under the isomorphism $\alpha: R\to S$ defined by
$u_j\mapsto y_j$ for $1\le j\le n$, we have
$\alpha(R^k)=S^k$ for all $k\ge 1$. Also
$\alpha(\mathcal I^k_\lambda)\equiv I^k_\lambda\mod S^{k-1}$ for all
$k\ge 0$.  Therefore $\alpha$ induces an isomorphism of abelian groups
\[{S_k}/{I_{\lambda,k}}\cong {S^k}/{(I^k_{\lambda}+S^{k-1})}\cong {R^k}/{(\mathcal
    I^k_\lambda+R^{k-1})}.\] Since by Theorem \ref{tanisaki}, for
every $k$ we have
\[{S_k}/{I_{\lambda,k}}\cong H^k(\mathcal F_\lambda;\mathbb Z),\] 
we conclude that, as abelian groups, 
\[{R}/{\mathcal I_{\lambda}}\cong H^*(\mathcal F_\lambda;\mathbb Z).\]  This
completes the proof.
\end{proof}
\noindent {\bf Acknowledgments:} We thank Professor Hiraku Abe for informing 
us about the work of Spaltenstein \cite{spal} and also for sending us
a copy of it. We thank the referee for a careful reading of the
manuscript and for valuable comments and suggestions.

\end{document}